\documentclass[12pt]{article}
\usepackage[latin1]{inputenc}
\usepackage{amssymb}
\usepackage{amsmath, amsthm}
\usepackage{amsfonts}
\numberwithin{equation}{section}

\usepackage[english]{babel}
\newtheorem{thm}{Theorem}[section]
\newtheorem{prop}[thm]{Proposition}

\newtheorem{lem}[thm]{Lemma}
\newtheorem{rem}{Remark}
\theoremstyle{definition}
\newtheorem{exa}{Example}
\newtheorem{defi}{Definition}[section]

\title{Common fixed points for pairs of mappings with variable contractive parameters}

\author{J.R. Morales and E.M. Rojas}

\date{}
\begin{document}
\maketitle
%\tikz[overlay,remember picture]
%{
%    \node at ($(current page.west)+(1.5,0)$) [rotate=90] {\Huge\textcolor{gray}{\textsf{Nueva Versi\'on \today\ a las \currenttime}}};
%    }

\begin{abstract}
In  this paper, we establish some common fixed point results for a new class of pair of contractions mappings having
functions as contractive parameters, and satisfying certain commutative properties. 
\end{abstract}

\section{Introduction and preliminaries}

The metric fixed point theory has a vast literature, the Banach-Caccioppoli's contraction principle is one of the most outstanding result in this theory. Since its appearance,
several generalizations of this result have been appeared in the literature.  In 1976, Jungck \cite{J76} generalize this principle by considering two commutating mappings
 and proved a common fixed point theorem for these mappings. Afterwards,  the commutative property of the mappings
assumed by Jungck has been relaxed by introducing ``weak'' alternative notions as \textit{weakly commutativity}, (\textit{non}-) \textit{compatible} and
 \textit{R-weakly commutating} among others, which allowed to generalize different classes of well-known contractive type of mappings. 

Here, we are going  to establish existence and uniqueness results of common fixed points  for a pair of contractive-type of mappings whose contractive parameters are non-constants 
and its contractive's inequality  is controlled by altering distance functions. To attain our goals,  we will assume that the mappings under consideration
are  occasionally weakly compatible. Also, alternatively, we will assume that the pair of mappings satisfies the so-called property (E.A.).
 
We would like to recall that  in 1984, M.S. Khan et al \cite{KSS84} introduced the notion of altering distance functions. Since then,
 it has been used to solve several problems in the metric fixed point theory (see, e.g., \cite{DuCh08,DuBiChDa09,MR12,MR12a,Na03,NSK11,SaBa98}). 
 
 \begin{defi}
A function $\psi: \mathbb{R}_{+}\longrightarrow \mathbb{R}_{+}:=[0,+\infty)$ is called an altering distance function if the following properties are satisfied:
\begin{enumerate}
\item[$(\Psi_1)$] $\psi(t)=0$ iff $t=0$.
\item[$(\Psi_2)$] $\psi$ is monotonically non-decreasing.
\item[$(\Psi_3)$] $\psi$ is continuous.
\end{enumerate}
By $\Psi$ we are going to denote the set of all the altering distance functions.
\end{defi}

In 2002 A. Branciari \cite{Br02} extends the Banach-Caccioppoli's theorem by using some Lebesgue integrable functions. Since then,
several well known fixed point criteria for contractive type of mappings have been generalized in this way. See e.g.,  \cite{Al06,Ay,CRSSS12,LLMC,MR12,MR12a,Rh03,SV12,SK11} and a lot of references therein.

By $\Phi$ is denoted the set of all mappings $\phi: \mathbb{R}_{+}\longrightarrow \mathbb{R}_{+}$ satisfying the following conditions:
\begin{enumerate}
\item[$(\Phi_1)$] $\phi$ is a Lebesgue integrable mapping which is summable on each compact subset of $\mathbb{R}_+$.
\item[$(\Phi_2)$] $\phi$ is non negative.
\item[$(\Phi_3)$] For each $\varepsilon>0$, $\int_0^\varepsilon \phi(t)dt>0$. 
\end{enumerate}
A relation between these two classes of functions $\Psi$ and $\Phi$ is given in the following result (\cite{AK09,PM09}):
\begin{lem}\label{lem:composition}
For each $\varphi\in\Phi$, the function $\psi_0:\mathbb{R}_+\longrightarrow \mathbb{R}_+$ defined by $\psi_0(s)=\int_0^s\varphi(t)dt$, $s\in \mathbb{R}_+$ is such that $\psi_0\in \Psi$.  
\end{lem}

In order to establish our results the following notions will be needed: A pair of selfmappings $(S,T)$ on a metric space $(M,d)$ is said \textit{compatible} \cite{J86} if and only if $\lim_{n\to\infty}d(TSx_n,STx_n)=0$, whenever
$(x_n)_n\subset M$ is such that 
\begin{equation*}%\label{compatible}
\lim_{n\to\infty}Sx_n=\lim_{n\to\infty}Tx_n=t
\end{equation*}
for some $t\in M$. A pair of mappings $(S,T)$
is said to be \textit{noncompatible} \cite{AM02} if there exits at least one sequence $(x_n)_n\subset M$ such that $\lim_{n\to\infty}Sx_n=\lim_{n\to\infty}Tx_n=t$
for some $t\in M$, but $\lim_{n\to\infty}d(STx_n,TSx_n)$ is either nonzero or non-existent.  A pair of selfmappings $(S,T)$ is said to satisfy the \textit{property} (E.A.), \cite{AM02}, if there exists a sequence $(x_n)_n\subset M$ such that 
\begin{equation*}%\label{property (E.A.)}
\lim_{n\to\infty}Sx_n=\lim_{n\to \infty} Tx_n=t,
\end{equation*} 
for some $t\in M$. 

 A point $x\in M$ is called a \textit{coincidence point} (CP) of $S$ and $T$ if $Sx=Tx$. The set of coincidence points of $S$ and $T$
will be denoted by $C(S,T)$. If $x\in C(S,T)$, then $w=Sx=Tx$ is called a \textit{point of coincidence} (POC) of $S$ and $T$. 

Finally, a pair of mappings $(S,T)$ is said to be \textit{occasionally weakly compatible} (OWC), \cite{TS08}, if  there exists $x\in C(S,T)$
such that $STx=TSx$. 

\begin{rem}\label{rem:CFP}
 In 2006, G. Jungck and B.E. Rhoades \cite{JR06} proved that if a pair of OWC maps $(S,T)$ has a unique POC, then it has a unique common fixed point, \label{un.CFP} and
M. Abbas and B.E. Rhoades \cite{AR08} asserted that the property (E.A.) implies OWC.

G.U.R. Babu and G.N. Alemyehu in \cite{BA10} proved that every pair of noncompatible selfmaps on a metric space satisfies the property (E.A.), but its converse 
is not true. Also, they showed that the property (E.A.) and OWC are independent conditions.
\end{rem}

The following lemma due to G.U. Babu and P.P. Sailaja in \cite{BS11} will be useful in the sequel.

\begin{lem}\label{lem1.5}
Let $(M,d)$ be a metric space. Let $(x_n)$ be a sequence in $M$ such that
\begin{equation*}%\label{eq1.4}
\lim_{n\rightarrow \infty} d(x_n, x_{n+1})=0.
\end{equation*} 
If $(x_n)$ is not a Cauchy sequence in $M$, then there exist an $\varepsilon_0>0$ and sequences of integers positive $(m(k))$ and $(n(k))$ with
\begin{equation*}
m(k)>n(k)>k
\end{equation*}
 such that,
 \begin{equation*}
d(x_{m(k)},\, x_{n(k)})\geq \varepsilon_0,\quad d(x_{m(k)-1},\, x_{n(k)})<\varepsilon_0
 \end{equation*}
 and
\begin{enumerate}
\item[(i)] $\lim_{k\rightarrow \infty}d(x_{m(k)-1},\, x_{n(k)+1})=\varepsilon_0$,
\item[(ii)] $\lim_{k\rightarrow \infty}d(x_{m(k)},\, x_{n(k)})=\varepsilon_0$,
\item[(iii)] $\lim_{k\rightarrow \infty}d(x_{m(k)-1},\, x_{n(k)})=\varepsilon_0$.
\end{enumerate}
\end{lem}

%\begin{rem}\label{ref1.6}
%From Lemma \ref{lem1.5} is easy to get that
%\begin{equation*}
%\lim_{k\rightarrow \infty} d(x_{m(k)+1},\, x_{n(k)+1})=\varepsilon_0.
%\end{equation*}
%\end{rem}

\section{The class of pairs of mappings with non-constant contractive parameters}

In order to introduce the class of mappings which will be the focus of study of this paper, as in \cite{LLMC}, we are going to use the functions $\alpha,\beta,\gamma:\mathbb{R}_+\longrightarrow [0,1)$ which satisfy that $\alpha(t)+\beta(t)+\gamma(t)<1$, for all $t\in \mathbb{R}_+$,
 and
\begin{align}\label{functions alpha, beta}
\limsup_{s\to 0^+}\gamma(s)<&1 \nonumber\\
\limsup_{s\to t^+}\frac{\alpha(s)+\beta(s)}{1-\gamma(s)}<&1,\quad \forall t>0.
\end{align}  

Now, we introduce the following class of pair of contraction-type of mappings.

\begin{defi}
Let $(M,d)$ be a metric space and let $S,T:M\longrightarrow M$ be mappings. The pair $(S,T)$ is called a $\psi-(\alpha,\beta,\gamma)$-contraction pair if for all
$x,y\in M$ 
\begin{align}\label{eq:1.3}
\psi\left(d(Sx,Sy)\right)\leq& \alpha( d(Tx,Ty))\psi(d(Tx,Ty))+\beta(d(Tx,Ty))\psi(d(Sx,Tx))\nonumber\\
                                          &+\gamma(d(Tx,Ty))\psi(d(Sy,Ty))
\end{align}
where $\psi\in \Psi$ and $\alpha,\beta,\gamma:\mathbb{R}_+\longrightarrow[0,1)$ are functions satisfying the conditions \eqref{functions alpha, beta}.
\end{defi}

\begin{rem}
In addition, we can consider a class of pair of mappings satisfying the following inequality contraction of integral type:
\begin{align}\label{int ineq}
\int_0^{\psi(d(Sx,Sy))}\varphi(t)dt\leq&\alpha(d(Tx,Ty))\int_0^{\psi(d(Tx,Ty))}\varphi(t)dt+\beta(d(Tx,Ty))\nonumber\\
&\times\int_0^{d(Sx,Tx)}\varphi(t)dt+\gamma(d(Tx,Ty))\int_0^{\psi(d(Sy,Ty))}\varphi(t)dt
\end{align}
for all $x,y\in M$, where $\psi\in\Psi$, $\varphi\in \Phi$ and $\alpha,\beta,\gamma$ are functions satisfying \eqref{functions alpha, beta}.
Notice that this class can be rewrite as 
\begin{align*}%\label{eq:1.3}
&\psi_0(\psi\left(d(Sx,Sy)\right))\leq \alpha( d(Tx,Ty))\psi_0(\psi(d(Tx,Ty)))\\
&+\beta(d(Tx,Ty))\psi_0(\psi(d(Sx,Tx)))+\gamma(d(Tx,Ty))\psi_0(\psi(d(Sy,Ty)))
\end{align*}
for all $x,y\in M$, where $\psi_0\in \Psi$ is the function defined in the Lemma \ref{lem:composition}. Thus, we have that all the conclusions 
given for $\psi-(\alpha,\beta,\gamma)$-contraction pairs are valid for pair of mappings satisfying the inequality contraction \eqref{int ineq}.
\end{rem}

%\begin{exa}\label{exa:1}
%Let $(M,d)$ be a metric space. If we consider $S\equiv a$, a constant map, and $T$ any selfmapping on $M$, we can check that the pair
%$(S,T)$ is a  $\psi-(\alpha,\beta,m)$-contraction pair for all functions $\alpha,\beta:\mathbb{R}_+\longrightarrow [0,1)$ such that $\alpha(t)+\beta(t)<1$,
% for all $t\in \mathbb{R}_+$ and satisfying \eqref{functions alpha, beta}.
%\end{exa}

\begin{prop}\label{lem:1.11}
Let $S$ and $T$ be two selfmaps on a metric space $(M,d)$. Let us assume that the pair $(S,T)$ is a $\psi-(\alpha,\beta,\gamma)$-contraction pair. If $S$
and $T$ have a POC in $M$ then it is unique.
 \end{prop}
 \begin{proof}
 Let $z\in M$ be a POC of the pair $(S,T)$. Then there exits $x\in M$ such that $Sx=Tx=z$. Suppose that for some $v\in M$, $Sv=Tv=w$ with $z\neq w$.
 Then, by \eqref{eq:1.3} we have
 \begin{align*}%\label{eq:1.10}
 &\psi(d(z,w))=\psi(d(Sx,Sv)) \leq \alpha(d(Tx,Tv))\psi(d(Tx,Tv))\\
 & +\beta(d(Tx,Tv))\psi(d(Sx,Tx))+\gamma(d(Tx,Tv))\psi(d(Sv,Tv)).
\end{align*}  
It follows that,
\begin{align*}%\label{eq:1.11}
\psi(d(z,w))\leq \alpha(d(z,w))\psi(d(z,w))
+\beta(d(z,w))\psi(d(z,z))+\gamma(d(z,w))\psi(w,w).
\end{align*}
Thus, we get
\begin{equation*}
\psi(d(z,w))\leq \alpha(d(z,w))\psi(d(z,w))<\psi(d(z,w).
\end{equation*}
since $\psi\in \Psi$, then we have that $d(z,w)<d(z,w)$ which is a contradiction, therefore $z=w$.
 \end{proof}

\begin{prop}\label{prop:1.2}
Let $(M,d)$ be a metric space and let $S,T:M\longrightarrow M$ be mappings with $S(M)\subset T(M)$. If the pair $(S,T)$ is a 
$\psi-(\alpha,\beta,\gamma)$-contraction pair, then for any $x_0\in M$, the sequence $(y_n)$ defined by
\begin{equation*}%\label{eq:1.13}
y_n=Sx_n=Tx_{n+1},\quad n=0,1,\dots
\end{equation*}
satisfies:
\begin{enumerate}
 \item[(1)] $\lim_{n\to\infty}d(y_n,y_{n+1})=0$.
 \item[(2)] $(y_n)\subset M$ is a Cauchy sequence in $M$. 
 \end{enumerate} 
\end{prop}
\begin{proof}
To prove (1), let $x_0\in M$ be an arbitrary point. Since $S(M)\subset T(M)$, then there exists $x_1\in M$ such that $Sx_0=Tx_1$. By continuing this process inductively
we obtain a sequence $(x_n)$ in $M$ such that
\begin{equation*}
y_n=Sx_n=Tx_{n+1}.
\end{equation*}  
Now, we have
\begin{align*}%\label{eq:1.14}
&\psi(d(Tx_{n+1},Tx_{n+2}))=\psi(d(Sx_{n},Sx_{n+1}))\leq \nonumber\\
&\alpha(d(Tx_{n},Tx_{n+1}))\psi(d(Tx_{n},Tx_{n+1}))+\beta(d(Tx_{n},Tx_{n+1}))\psi(d(Sx_{n},Tx_{n}))\\
&+\gamma(d(Tx_n,Tx_{n+1})\psi(d(Sx_{n+1},Tx_{n+1})).
\end{align*}
It follows that
\begin{align*}
&\psi(d(Tx_{n+1},Tx_{n+2}))\leq \alpha(d(Tx_n,Tx_{n+1}))\psi(d(Tx_n,Tx_{n+1}))\\
&+\beta(d(Tx_n,Tx_{n+1}))\psi(d(Tx_{n+1},Tx_n))+\gamma(d(Tx_{n},Tx_{n+1})\psi(d(Tx_{n+2},Tx_{n+1})).
\end{align*}
Therefore, we obtain
\begin{align*}%\label{eq:1.15}
\psi(d(Tx_{n+1},Tx_{n+2}))\leq \frac{\alpha(d(Tx_{n},Tx_{n+1}))+\beta(d(Tx_n,Tx_{n+1}))}{1-\gamma(d(Tx_{n},Tx_{n+1}))}\psi(d(Tx_{n},Tx_{n+1})),
\end{align*}
from which, together with \eqref{functions alpha, beta}, we conclude that
\begin{equation*}
\psi(d(Tx_{n+1},Tx_{n+2}))<\psi(d(Tx_n,Tx_{n+1})).
\end{equation*}
  Since $\psi\in \Psi$, then $(d(Tx_{n},Tx_{n+1}))$ is a monotone decreasing sequence of non negative real numbers which converges to $a\geq 0$. Thus,
\begin{equation*}
\lim_{n\to\infty}d(Tx_{n},Tx_{n+1})=a\geq 0.
\end{equation*}  
We want to prove that $a\equiv0$. We are going to assume that $a>0$. Using the continuity of $\psi$ and conditions \eqref{functions alpha, beta} we have
\begin{align*}
&0<\psi(a)=\limsup_{n\to\infty}\psi(d(Tx_n,Tx_{n+1}))\\
    &\leq\limsup_{n\to\infty}\left(\frac{\alpha(d(Tx_{n-1},Tx_n)+\beta(Tx_{n-1},Tx_n))}{1-\gamma(d(Tx_{n-1},Tx_n))}\right)\psi(d(Tx_{n-1},Tx_n))<\psi(a)
\end{align*}
which is  a contradiction. So, $a=0$ and therefore,
\begin{equation*}
\lim_{n\to\infty}\psi(d(Tx_n,Tx_{n+1}))=0,
\end{equation*}
thus,
\begin{equation*}%\label{eq:1.17}
\lim_{n\to\infty}d(y_n,y_{n+1})=\lim_{n\to\infty}d(Sx_{n-1},Sx_{n})=\lim_{n\to\infty}d(Tx_{n},Tx_{n+1})=0.
\end{equation*}
To prove (2), we are going to suppose that $(y_n)\subset T(M)$ is not a Cauchy sequence. Then, from Lemma \ref{lem1.5} there exists $\varepsilon>0$ and 
sequences $(m(k))$ and $(n(k))$ with $m(k)\geq n(k)>k$ such that 
\begin{equation*}
\lim_{k\to\infty}d(Tx_{m(k)},Tx_{n(k)})=\varepsilon
\end{equation*}
and
\begin{equation*}
\lim_{k\to\infty}d(Tx_{m(k)-1},Tx_{n(k)-1})=\varepsilon.
\end{equation*}

In this way we have
\begin{align*}%\label{eq:1.18}
0<\psi(\varepsilon)=&\limsup_{k\to\infty}\psi(d(Tx_{m(k)},Tx_{n(k)}))= \limsup_{k\to\infty}\psi(d(Sx_{m(k)-1},Sx_{n(k)-1}))\nonumber\\
                                  \leq& \limsup_{k\to\infty}\alpha(d(Tx_{m(k)-1},Tx_{n(k)-1}))\psi(d(Tx_{m(k)-1},Tx_{n(k)-1}))\nonumber\\
                                  &+\limsup_{k\to\infty}\beta(d(Tx_{m(k)-1},Tx_{n(k)-1}))\psi(d(Sx_{m(k)-1},Tx_{m(k)-1}))\\
                                  &+\limsup_{k\to\infty}\gamma(d(Tx_{m(k)-1},Tx_{n(k)-1}))\psi(d(Sx_{n(k)-1},Tx_{n(k)-1}))\\
                                 \leq & \limsup_{s\to\varepsilon}\alpha(s)\limsup_{k\to\infty}\psi(d(Tx_{m(k)-1},Tx_{n(k)-1}))\nonumber\\
                                  &+\limsup_{s\to\varepsilon}\beta(s)\limsup_{k\to\infty}\psi(d(Tx_{m(k)},Tx_{m(k)-1}))\\
                                  &+\limsup_{s\to\varepsilon}\gamma(s)\limsup_{k\to\infty}\psi(d(Tx_{n(k)},Tx_{n(k)-1}))<\psi(\varepsilon),
\end{align*}
which is a contradiction, hence $(Tx_n)$ is a Cauchy sequence in $M$.
\end{proof}

\section{On the existence and uniqueness of common fixed points}

In this section we prove our main results concerning to the existence and uniqueness of common fixed points for a  $\psi-(\alpha,\beta,\gamma)$-contraction pair of mappings without continuity requirement. 
In comparison with classical results in this theory, the commutativity property in this case is reduced to the existence of points of
coincidence and the completeness of the space is reduced to natural conditions.

\begin{thm}\label{thm:2.1}
Let $S$ and $T$ be selfmaps on a metric space $(M,d)$ such that
\begin{enumerate}
\item[(i.)] $S(M)\subset T(M)$.
\item[(ii.)] $T(M)\subset M$ is a complete subspace of $M$.
\item[(iii.)] The pair $(S,T)$ is a $\psi-(\alpha,\beta,\gamma)$-contraction pair. 
\end{enumerate}
Then,
\begin{enumerate}
\item[(1)] The pair $(S,T)$ has a unique POC.
\item[(2)] If the pair $(S,T)$ is OWC, then $S$ and $T$ have a unique common fixed point. 
\end{enumerate}
\end{thm}
\begin{proof}
Let $y_n=Sx_n=Tx_{n+1}$, $n=0,1,\dots,$ be the Cauchy sequence defined in Proposition \ref{prop:1.2} which, as was proved, satisfies that
 $(y_n)=(Tx_{n+1})\subset T(M)$.  Since $T(M)\subset M$ is a complete subspace of $M$, then there exists $z\in T(M)$ such that
\begin{equation*}%\label{eq:2.1}
\lim_{n\to\infty}y_n=\lim_{n\to\infty}Sx_n=\lim_{n\to\infty}Tx_{n+1}=z,
\end{equation*}
and thus we can find $u\in M$ such that $Tu=z$. Now, we are going to assume that $Su\neq z$. Then,
\begin{align*}
&\psi(d(Sx_{n+1},Su))\leq \alpha(d(Tx_{n+1},Tu))\psi(d(Tx_{n+1},Tu))\\
&+\beta(d(Tx_{n+1},Tu))\psi(Sx_{n+1},Tx_{n+1})+\gamma(Tx_{n+1},Tu)\psi(Su,Tu).
\end{align*}
Letting $n\to\infty$, we obtain  
\begin{align*}
&\psi(d(Su,z))\leq \limsup_{n\to\infty}\alpha(d(Tx_{n+1},Tu))\psi(d(z,Tu))\\
&+ \limsup_{n\to\infty}\beta(d(Tx_{n+1},Tu))\psi(d(z,z))+ \limsup_{n\to\infty}\gamma(d(Tx_{n+1},Tu))\psi(d(Su,Tu))\\
&<\psi(d(Su,z))
\end{align*}
which is a contradiction, therefore $Su=z$, hence $z$ is a POC of $S$ and $T$. 
 From the Proposition \ref{lem:1.11} we conclude that $z$ is the unique POC. 

On the other hand, since the pair $(S,T)$ is OWC, then it has a unique common fixed point (see Remark \ref{rem:CFP}).  
\end{proof}

\begin{thm}\label{thm:fin}
Let $(M,d)$ be a metric space and $S,T:M\longrightarrow M$ mappings satisfying the property (E.A.). Let us suppose that the pair $(S,T)$
is a $\psi-(\alpha,\beta,\gamma)$-contraction pair. If $T(M)\subset M$ is closed, then $S$ and $T$ have a unique common fixed point.
\end{thm}
\begin{proof}
Since the pair $(S,T)$ satisfies the property (E.A.), then there exists a sequence $(x_n)\subset M$ such that 
\begin{equation*}
\lim_{n\to\infty}Sx_n=\lim_{n\to\infty} Tx_n=z
\end{equation*}
  for some $z\in M$. Since $T(M)$ is closed, then $z\in T(M)$ and $z=Tu$ for some $u\in M$. As in the proof of the Theorem \ref{thm:2.1}, we
  can prove that $z=Tu=Su$ and that $z$ is the unique POC of $S$ and $T$. 
  Finally, since the pair $(S,T)$ is OWC, then $z$ is the unique common fixed point.
\end{proof}

%\begin{exa}
%Let us consider the mappings of the Example \ref{exa:1}. Furthermore, let us assume that $\{a\}\subset T(M)$.  Then, $a$ is the unique
%POC of the pair $(S,T)$. Notice that $(S,T)$ is WOC if and only if $a$ is a fixed point of $T$. Therefore, $a$ is
%the unique common fixed point of the pair $(S,T)$. We would like to point out that $(S,T)$ satisfies the property (E.A.), for the constant sequence $x_n\equiv a$.    
%\end{exa}

\begin{rem}
Since two noncompatible selfmappings on a metric space $(M,d)$ satisfy the property (E.A.), then the conclusion of the Theorem \ref{thm:fin} remains valid
if we consider $S$ and $T$ noncompatible selfmappings. 
\end{rem}

\section{Conclusions and examples}

Notice that our results extend several classes of well known contractive type of mappings, including various classes of contractive mappings of the integral type. Even more, the mappings $T$ and $S$ considered
here are not necessarily  continuous, so in this way our results are more general compared with other results in this line of research.  

Next, we are going to show some examples in support of our results.

\begin{exa}
Let be $M=[0,1]$ equipped with the euclidean metric. We consider the following mappings: $S,T:M\longrightarrow M$ defined by $Sx=\frac{x}{16}$ and $Tx=\frac{x}{2}$ for all
$x\in M$. Let $\alpha,\beta,\gamma:\mathbb{R}_+\longrightarrow[0,1)$ defined by $\alpha(t)=\beta(t)=\frac{1}{2}$ and $\gamma(t)=\frac{1}{4}$ for all $x\in \mathbb{R}_+$ and $\psi:\mathbb{R}_+\longrightarrow\mathbb{R}_+$ given by the formula $\psi(t)=t^2$, $t\in\mathbb{R}_+$.

Notice that $\psi\in \Psi$ and the functions $\alpha,\beta,\gamma$ satisfy the conditions \eqref{functions alpha, beta}, also note that $S(M)\subset T(M)$ and $T(M)$ is a complete subspace of $M$. Moreover, 
is not difficult to show that the pair $(S,T)$ is a $\psi-(\alpha,\beta,\gamma)$-contraction pair. Besides $C(S,T)=\{0\}$ and $ST0=TS0=0$, which mean that $(S,T)$ is OWC. Then, the Theorem \ref{thm:2.1} guarantee that $w=0$ is the unique common fixed point of $S$ and $T$.
\end{exa}

\begin{exa}
As in the example before, $M=[0,1]$ with the usual metric. We define the selfmaps $S,T$ on $M$ by
\begin{equation*}
Sx=\begin{cases}
0, & \mbox{if } 0\leq x\leq\frac{1}{2} \\ 
\frac{1}{16}, & \mbox{if } \frac{1}{2}< x\leq1
 \end{cases}
\end{equation*} 
and $Tx=\frac{x}{2}$ for all $x\in M$. Let $\alpha,\beta,\gamma:\mathbb{R}_+\longrightarrow[0,1)$ defined as follow:
\begin{equation*}
\alpha(t)=\frac{1}{8},\quad \beta(t)=\gamma(t)=\frac{1}{4},\quad \mbox{for all $t\in\mathbb{R}_+$}.
\end{equation*}
Let  $\psi:\mathbb{R}_+\longrightarrow\mathbb{R}_+$ defined by $\psi(t)=t^2$, $t\in\mathbb{R}_+$. Then, the pair $(S,T)$ is a $\psi-(\alpha,\beta,\gamma)$-contraction pair
satisfying the hypotheses of Theorem \ref{thm:2.1}, thus $w=0$ is the unique POC and moreover the unique common fixed point of $S$ and $T$. 
\end{exa}

\begin{exa}
Let $M=[1/2,1]$ with the usual metric on $\mathbb{R}$. In this case we consider the mappings $S,T:M\longrightarrow M$ defined by
\begin{equation*}
Sx=\begin{cases}
\frac{1}{2}, & \mbox{if } \frac{1}{2}\leq x<\frac{2}{3} \\ 
1-\frac{1}{2}x, & \mbox{if } \frac{2}{3}\leq x\leq1
 \end{cases},\quad
  Tx=\begin{cases}
1, & \mbox{if } \frac{1}{2}\leq x<\frac{2}{3} \\ 
x, & \mbox{if } \frac{2}{3}\leq x\leq1
 \end{cases}.
\end{equation*} 
Let $\psi:\mathbb{R}_+\longrightarrow\mathbb{R}_+$ given by the formula $\psi(t)=t^2$, $t\in\mathbb{R}_+$. $\alpha,\beta,\gamma:\mathbb{R}_+\longrightarrow [0,1)$
given by $\alpha(t)=\beta(t)=\frac{1}{4}$, $\gamma(t)=\frac{1}{8}$, $t\in\mathbb{R}_+$. Notice that $C(S,T)=\{2/3\}$, and $ST(2/3)=TS(2/3)$. Moreover, $w=2/3$ is the unique POC of $S$ and $T$, thus the pair $(S,T)$ is OWC. On the other hand, by considering the sequence $x_n=\frac{2}{3}+\frac{1}{n}$, $n\geq 4$ in $M$, it is clear that the pair $(S,T)$
satisfies the property (E.A.). Finally, it is easy to show that in fact the mappings $S,T,\psi,\alpha,\beta$ and $\gamma$ satisfy all the hypotheses of Theorem \ref{thm:fin}, so $w=2/3$ is the unique common fixed point of $S$ and $T$.

\end{exa}

\noindent
J.R. Morales, Departamento de Matem\'aticas, Universidad de Los Andes, M\'erida 5101, Venezuela, \texttt{moralesj@ula.ve}\\

\noindent
E.M. Rojas, Departamento de Matem\'aticas, Pontificia Universidad Javeriana, Bogot\'a, Colombia, \texttt{edixon.rojas@javeriana.edu.co}

\end{document}